\documentclass[10pt,conference, twocolumns]{IEEEtran}

\usepackage{times}
\usepackage{cite} 
\usepackage[english]{babel}
\usepackage{epsfig} 
\usepackage{amsfonts}
\usepackage{bm}
\usepackage{amsbsy}
\usepackage{amsmath} 
\usepackage{graphicx}
\usepackage{balance} 
\usepackage{latexsym}
\usepackage{amssymb}
\usepackage{txfonts} 
\usepackage{wasysym}
\usepackage{mathbbol}
\usepackage{multirow} 
\usepackage{epstopdf} 
\usepackage[ruled,vlined]{algorithm2e} 
\usepackage{color}

\newtheorem{theorem}{Theorem}
\newtheorem{corollary}{Corollary}

\newtheorem{remark}{Remark}
\usepackage{authblk}
  \newcommand{\prob}{\mathrm{\textit{P}}}

\DeclareMathOperator{\sign}{\textit{Sign}}
\allowdisplaybreaks
\title{Minimum Throughput Maximization in LoRa Networks Powered by Ambient Energy Harvesting}
\author{Fatma Benkhelifa\**, Zhijin Qin\***, Julie McCann\**\\
\small \** Imperial College London, London, UK\\
\small \*** Queen Mary University of London, London, UK\\
\small \{f.benkhelifa,j.mccann\}@imperial.ac.uk, z.qin@qmul.ac.uk
}

\begin{document}
\maketitle

\begin{abstract}
    {\ In this paper, we investigate the uplink transmissions in low-power wide-area networks (LPWAN) where the users are self-powered by the energy harvested from the ambient environment. Demonstrating their potential in supporting diverse Internet-of-Things (IoT) applications, we focus on long range (LoRa) networks where the LoRa users are using the harvested energy to transmit data to a gateway via different spreading codes. Precisely, we study the throughput fairness optimization problem for LoRa users by jointly optimizing the spreading factor (SF) assignment, energy harvesting (EH) time duration, and the transmit power of LoRa users. First, through examination of the various permutations of collisions among users, we derive a general expression of the packet collision time between LoRa users, which depends on the SFs and EH duration requirements. Then, after reviewing prior SF allocation work, we develop two types of algorithms that either assure fair SF assignment indeed purposefully 'unfair' allocation schemes for the LoRa users. Our results unearth three new findings. Firstly, we demonstrate that, to maximize the minimum rate, the unfair SF allocation algorithm outperforms the other approaches. Secondly, considering the derived expression of packet collision between simultaneous users, we are now able to improve the performance of the minimum rate of LoRa users and show that it is protected from inter-SF interference which occurs between users with different SFs. That is, imperfect SF orthogonality has no impact on minimum rate performance. Finally, we have observed that co-SF interference is the main limitation in the throughput performance, and not the energy scarcity.
    \par}
\end{abstract}

\section{Introduction}
    {\ Internet of Things (IoT), sensor networks and cyber-physical systems have gained a lot of recent interest as they are used in a wide range of applications to enable their automation, resource optimization, resilience to change and failure and make them more sustainable. Examples include Smart Cities, Smart Farming/Agriculture, Health Care, Public Safety, etc~\cite{IOT1}. Many such systems involve the use of relatively large numbers of deployed devices with limited resources that need to deliver reliable data to potentially critical applications. Many of these systems span multikilometers, e.g. smart water networks, precision farms, etc., which have motivated the development of communication technologies that have low power operation. To achieve this, they operate at relatively low data rates. These low power wide area networks (LPWANs) are gaining serious interest across the world with countries such as the Netherlands deploying a LPWAN country-wide~\cite{lpwa}. Many standard technologies were competing to model LPWAN such as SigFox, Weightless, Narrowband IoT (NB-IoT) and Long Range (LoRa)~\cite{loraweb}. Among the existing standards for LPWAN, LoRa has captured a lot of research and industrial attention by its strength lying in its ability to cover large geographical distances and in Chirp Spread Spectrum (CSS) modulation that assigns different spreading factors (SFs) to users making it resilient to external interference~\cite{loraweb}. So far, most research interest in LoRa networks was more about the study of scalability, coverage, and reliability considering the co-SF and/or inter-SF interference.
    \par}
    {\ For example in~\cite{canlorascale}, the uplink coverage probability of a single LoRa gateway has been derived using stochastic geometry and two-link outage conditions were solved showing that the scalability of LoRa is related to co-spreading factor interference. In~\cite{imperfectSF}, a theoretical analysis of the achievable LoRa throughput has been analyzed taking into consideration the co-SF interference between users using the same SF as well as the inter-SF interference between users using different SFs. In~\cite{SF_collision_distance}, the average system packet success probability (PSP) has been analyzed using stochastic geometry for a LoRa system using the unslotted ALOHA random access protocol. Also, an adaptive SF allocation algorithm has been proposed which maximizes the average system PSP. 
    In~\cite{Zhijin}, the resource allocation of uplink transmissions in LoRa networks has been investigated where the minimum throughput rate of LoRa users has been optimized through the joint optimization of the channel assignment and the power allocation using a many-to-one matching game. Within each channel, the transmit powers for the LoRa users have been optimized. 
    In~\cite{fair_alloc}, a fair adaptive data rate allocation has been proposed to achieve a fair collision probability among all used data rates and a transmission power control algorithm has been developed to exhibit data rate fairness among the users. 
    \par}
{\  Nevertheless, LoRa resilience is limited if the devices are powered by finite energy (battery-based) sources, which also limit where such devices can be deployed as the cost of battery replacement is higher for devices positioned in difficult to access or dangerous environments. 
Thus, energy efficiency has been addressed in various studies~\cite{enereff,Zhijin2} in order to extend the battery lifetime of sensor devices. However, improving energy efficiency is not sufficient in itself. 
Ambient energy harvesting (EH) is a viable alternative to ensure sustainable operation or at least elongated life-times. EH can be obtained from different sources such as solar energy, wind energy, electromagnetic energy, radio frequency (RF) energy~\cite{fundamsurveryEH}. The latter technique can be obtained from dedicated transmitters (that exist specifically to provide energy) or ambient energy (transmitters that exist in the environment already e.g. WiFi). 
    Nonetheless, few research works have analyzed energy harvesting in LoRa networks. 
    In~\cite{ehlora}, a battery-less LoRa wireless sensor has been proposed that monitors road conditions powered by the vibration energy harvested by an electromagnetic energy harvester based on a Halbach configuration. 
    In~\cite{ehlora2}, a novel floating device has been proposed with a multi-source energy harvesting technique which harvests solar and thermoelectric energy and its extension, ~\cite{ehlora3}, that focus on power reduction when in listening. These works are relatively rudimentary as they do not carry out energy budgeting or optimization nor provide energy neutral guarantees. 
    \par}
    {\ In this paper, we study the resource allocation that maximizes the minimum rate of the LoRa users harvesting energy from an external source to enable operation and uplink transmissions. Thus, we propose two types of algorithms that allocate the SFs either fairly or unfairly between the users. We compare the proposed algorithms to other related algorithms in the literature. Then, we optimize the EH time and the power allocation between the LoRa users. In the simulation results, we harvest the energy from ambient RF signals and we compare our proposed solution to different baseline schemes. 
    The contributions of this paper are summarized as follows:
        \begin{itemize}
            \item To the best of our knowledge, this paper is the first to address resource allocation in LoRa networks using the EH capabilities to power its uplink transmissions. The EH source can be of any type as long as it is independent of the LoRa frequency band. 
                        \item In our analysis, we propose a general framework that accounts for the interference between users transmitting at the same time over nonorthogonal waveform codes. Thus, the multiuser interference between the users occurs either due to colliding users in the time or due to nonorthogonality between spreading coded waveforms.
            \item We explicitly express the packet collision time between users using the same SF or different SFs, and we show that its expression depends on the EH time and SFs. 
            
                        \item We have observed that the unfair SF allocation outperforms all the other SF allocation approaches in terms of the minimum rate of LoRa users.
            \item We have observed that the inter-SF interference (if considered) is neutralized thanks to the expression of the packet collision time for a specific value of the EH time.
            \item Finally, we have seen that co-SF interference is the main limitation of the system performance, and not really the energy scarcity.
        \end{itemize}
    \par}

\section{System Model and Problem Formulation}
    {\ Recall, we consider the uplink transmissions from LoRa users that are self-powered by external energy harvesting sources. 
    LoRa wide-area network (LoRaWAN) employs a star-of-stars topology where gateways relay the data transmissions between the end user devices and the server. We assume that we have $U$ LoRa users that are uniformly distributed in a circle of radius $R$ centred around a Gateway. The channel between the n'th LoRa user and the gateway is modelled as a Rayleigh fading channel with path loss as
			\begin{align}
				g_n &= h_n d_n^{-\alpha},
			\end{align}
	where $h_n$ is the small scale fading that is exponentially distributed with unit mean, $d_n$ is the distance between the n'th LoRa user and the gateway, and $\alpha$ is the path loss exponent.
    \par}
\subsection{Physical Layer of LoRa}
    {\ LoRaWAN operates in the sub-GHz frequency bands~\cite{loraweb}. In Europe, LoRaWAN uses the EU industrial, scientific, and medical (ISM) $868$ MHz frequency. For this band, there are eight physical layers: six with SF from $7$ to $12$ with bandwidth $125$ kHz and one with SF $7$ on $250$ kHz and one with Gaussian frequency-shift keying (GFSK) at $50$ kbps data rate. For the medium access control (MAC) layer, the end user devices access the channel using the pure ALOHA for transmitting their packets. 
    In this work, the bandwidth $\text{BW}$ is equal to $125$ kHz for which we have 8 channels and maximum 6 SFs per channel. We assume that each LoRa user has always access to one free channel. \footnote{The channel access selection is out of the scope of this work and might be considered in future works.} The $n$'th LoRa user transmits during the time on air given by
	    \begin{align}
		    T_{a,n} &=nb_n \times \frac{2^{\text{SF}_n}}{\text{BW}},
		\end{align}
	where $nb_n$ is the number of symbols, and $\text{SF}_n\in \mathcal{S}= \{7,\dots,12\}$ is the spreading factor. Note that $T_{a,n}$ belongs to the space $\mathcal{T}_a= \{t_{a,i}= d \frac{ 2^i nb_i}{d \text{BW}}, i=7,\dots, 12 \}$ with dimension $6$. 
	In addition, the European frequency regulations impose duty cycle restrictions for the $868$ MHz sub-bands, either $1\%$ or $10\%$~\cite{doppler}. 
    Subsequently, each LoRa user should stay silent $(1-d)\%$ of the packet duration once he transmits over one channel, where $d$ is the duty cycle chosen equal to $1\%$. Thus, the time off per channel for the $n$'th users is expressed as $T_{off,n}= \frac{1-d}{d}  T_{a,n}$.
    \par}
\subsection{Energy Harvesting at LoRa Users}
    {\ Each LoRa user is batteryless and is powered by harvested energy from external sources of energy. 
    The external energy harvesting source could be of any type, under the condition that it is not interfering with the band frequency of the LoRa users. For example, if it is from radio frequency (RF) energy harvesting, we assume that the energy is harvested from a band frequency other than the $868$ MHz. 
    We assume that the harvested energy per time unit is known and uncontrollable. Let $E_n$ be the harvested energy per time unit of the $n$'th user. It is independent of the time and it depends only on external conditions such as its location compared to the energy sources (channel gain, distance, etc.). We consider the "harvest-then-transmit" protocol: each LoRa user harvests first what it needs and transmits its data later. For the $n$'th user, the harvested energy during a harvesting time $\tau_{e,n}$ is given by
    					\begin{align}
    						E_{h,n} &= \tau_{e,n} E_n.
    					\end{align}
    Since the LoRa user is the one that decides when to transmit, we can consider that each user can harvest energy during a time less or equal to $(1-d)\%$ of the packet duration before performing a transmission. This adds a constraint on the harvesting time
				\begin{align}
					0 \leq \tau_{e,n} \leq T_{off,n} = \frac{1-d}{d}   T_{a,n}.
				\end{align}
	Hence, the available power at each user after harvesting is given by
		\begin{align}
			P_{h,n} &=  \frac{E_{h,n}}{T_{a,n}} = \frac{\tau_{e,n} E_n}{ T_{a,n} }.
		\end{align}
	The LoRa user transmits with a maximum transmit power $P_t$ which is most likely to be known and we assume that it is the same for all users. We assume also that all users cannot store the remaining energy after transmission.  

	The extra energy is lost. The storage of energy will be discussed in future works. Subsequently, the power allocated per user $p_n$ is constrained to the maximum transmit power and to the available harvested power $P_{h,n}$ after harvesting.
    \par}
\section{Packet Collision Time Between LoRa Users}
    {\ In this section, we examine the packet collision time between users using either the same or different SFs. This means that we consider the general case where inter-SF interference between users using different SFs is also possible, namely the imperfect SF orthogonality case. Let $col_{nm}$ be the packet collision time between the user $n$ and the user $m$. 
    We consider the case where the users start transmitting immediately after finishing harvesting their energy, which means that the start time of transmission is exactly the same as the end of EH time. In order to derive the expression of the collision time between two users $n$ and $m$, we consider two cases: either user $n$ spends more time to harvest or less than user $m$. 
    \par} 
    {\ First, if the user $n$ requires more time to harvest than the user $m$, i.e. $\tau_{e,n} \geq \tau_{e,m}$, the packet collision time will depend on how much each user spends during the packet transmission, as shown in Fig. \ref{colltimeEH}. If the user $m$, who requires less time to harvest, finishes his packet transmission before the user $n$ starts transmitting, then the two users will not collide. On the other hand, if the user $m$ finishes its packet transmission after the user $n$ finishes its packet transmission, then the two users will collide during the packet transmission time $T_{a,n}$ of the user $n$. Otherwise, if the user $m$ finishes transmitting his packet after the user $n$ finished harvesting and before he finishes his packet transmission, then the collision time will be $\tau_{e,m}+T_{a,m}-\tau_{e,n}$. 
    Second, if the user $n$ requires less time to harvest than the user $m$, i.e. $\tau_{e,n} < \tau_{e,m}$, the collision time is expressed by just exchanging the roles of the users $n$ and $m$ above.
    In summary, the collision time between user $n$ and user $m$ depends on the EH periods, as well as the time on air, and its expression is given in \eqref{cola} and \eqref{colb}, with $\mathcal{M}_{n,m}= \max\left(T_{a,n},T_{a,m}\right)$ and $\mathcal{N}_{n,m}= \min\left(T_{a,n},T_{a,m}\right)$.
		\begin{itemize}
			\item If $\sign\left(\tau_{e,n} - \tau_{e,m}\right)=\sign\left(T_{a,n}-T_{a,m}\right)$,  
			\small{	\begin{align}
					col_{n,m}  \label{cola}
					&=\begin{cases}
					 	0, &\hspace{-3mm}\mbox{ if }  \vert \tau_{e,n} - \tau_{e,m} \vert\geq   \mathcal{N}_{n,m},\\
						\mathcal{N}_{n,m}-\vert \tau_{e,n} - \tau_{e,m} \vert, &\hspace{-3mm}\mbox{ if }  \vert\tau_{e,n}  - \tau_{e,m} \vert <  \mathcal{N}_{n,m}.
					\end{cases}
				\end{align}}
				\item If $\sign\left(\tau_{e,n}-\tau_{e,m}\right)\neq\sign\left(T_{a,n}-T_{a,m}\right)$, 
		\small{	\begin{align}
					col_{n,m}  \label{colb}
					&=\begin{cases}
					 	0, &\hspace{-3mm}\mbox{ if }  \vert \tau_{e,n} - \tau_{e,m} \vert \geq  \mathcal{M}_{n,m},\\
						 \mathcal{M}_{n,m}-\vert \tau_{e,n} - \tau_{e,m} \vert, &\hspace{-3mm}\mbox{ if }  \vert T_{a,m}-T_{a,n} \vert\leq \vert \tau_{e,n}  - \tau_{e,m} \vert <  \mathcal{M}_{n,m},\\
						 \mathcal{N}_{n,m}, &\hspace{-3mm}\mbox{ if }   \vert \tau_{e,n} - \tau_{e,m}  \vert <\vert T_{a,m} -T_{a,n} \vert.
					\end{cases}
				\end{align}}
		\end{itemize} 
        \begin{figure}[t]
       			 \begin{center}
      				  \includegraphics[scale=0.415]{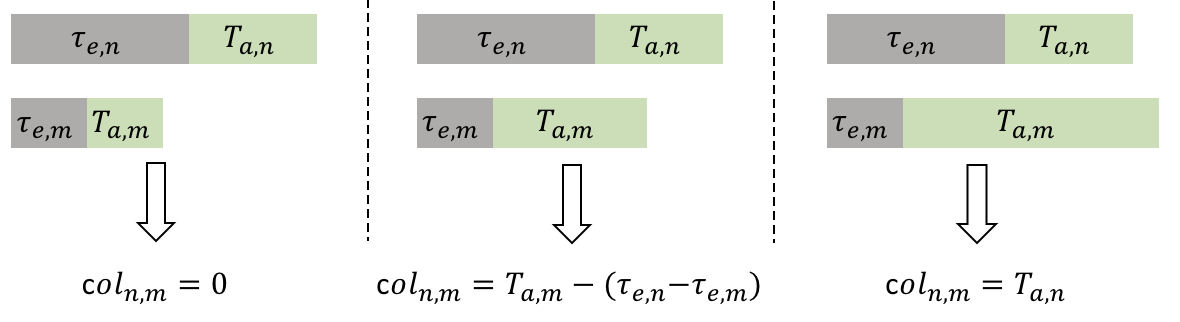}
      				  \caption{Collision time depending on EH time}
     			   \label{colltimeEH}
     			   \end{center}
    		 \end{figure} 
		 \begin{remark}
		 The collision time depending on the EH time is always less or equal to the collision time in the worst case scenario expressed as
		 \begin{align}
   			 col^{\text{worst case}}_{n,m}= \mathcal{N}_{n,m}.
   			 \label{colWC}
 		     \end{align}
		 This collision time can happen when all users finish transmitting at the same time which is the worst interference case. Here, the users with the lowest SFs will undergo the lowest interference (highest rate per unit of time) and the users with the highest SFs will undergo the highest interference (lowest rate per unit of time). 
	However, the users with the lowest SFs will have less time to transmit and harvest, while the users with the highest SFs will have more time to transmit and harvest.  
	 \end{remark}	
	 \begin{remark}\label{rk6}
	 If the EH time is expressed as $\tau_{e,n}=\frac{1-d}{d} T_{a,n}$, $\forall n$ and we satisfy the condition $\frac{t_{a,i}}{t_{a,j}}>1-d$, for $7\leq j < i \leq 12$, the collision time between the users simplifies to
	 \begin{align}
	 col_{n,m}=\begin{cases}
	 0, &\mbox{ if } \text{SF}_n\neq \text{SF}_m,\\
	 T_{a,n}, &\mbox{ if } \text{SF}_n=\text{SF}_m.
	 \end{cases}
	     \end{align}
	 This expression means that, even though we have imperfect orthogonality between SFs, the users using different SFs do not collide over the time domain.
	 \end{remark}
	\par}
    {\ Considering the packet collision time between users, the expression of signal-to-interference plus noise ratio (SINR) of user n is given by
	    \begin{align}
        	\gamma_n^{col} &= \frac{p_n g_n }{ \sum_{m \neq n } \frac{col_{n,m}}{T_{a,m}} \rho_{m,n}  p_m g_m + \sigma^2},
	    \end{align}
	where $\sigma^2$ is the variance of the additive white Gaussian noise (AWGN) at the LoRa gateway, and $\rho_{m,n}$ is the correlation factor between the coded waveforms for user m and user n. Note that $\rho_{m,n}=1$ if $m=n$ and $0\le\rho_{m,n}<1$ if $m \ne n$. Moreover, for the users sharing the same SF over the same channel, the correlation factor will be much higher than that of users sharing different SFs. 
    \par}
\section{Max-Min Throughput Optimization for Uplink Transmissions}
    {\ In this section, we propose to maximize the minimum transmission uplink rate of all the LoRa users while optimizing the SF allocation, the EH time assignment, and the power allocation between users.
        \begin{subequations}
	    	\begin{align}
	        	\underset{\tau_{e,n}, p_n, T_{a,n}}{\max}&\underset{n \in U}{\min} \hspace{3mm} R_n = \log\left(1 + \gamma_n^{col} \right),\\
		        \text{s.t. }
	        	& \text{ C$_{1}$: } 0 \leq p_n \leq P_t,  \\
	        	& \text{ C$_{2}$: }0 \leq p_n \leq \frac{\tau_{e,n} E_n}{T_{a,n}}, \\
	        	& \text{ C$_{3}$: } 0 \leq \tau_{e,n} \leq \frac{1-d}{d}T_{a,n}, \\
	        	& \text{ C$_{4}$: } T_{a,n} \in \mathcal{T}_a, 
		    \end{align}\label{mainProb}
    	\end{subequations}	
    The constraint C$_1$ is due to the maximum transmit power constraint at each LoRa user, the constraint C$_2$ is because each user cannot use a power more than the available harvested power, the constraint C$_3$ is due to the duty cycle restriction, and the constraint C$_4$ is due to the different SF assignment at each user. 
    This optimization problem \eqref{mainProb} is non-convex since it is a mixed-integer programming problem and the objective function is non concave due to the interference from colliding users. So, no computationally efficient method can be proposed without involving exhaustive search which has at least a complexity of the order of $U^{\vert \mathcal{T}_a\vert}N_{\epsilon}^2$, where $N_{\epsilon}$ is the complexity of the one-dimensional search method. To simplify the analysis, we propose to first decouple the problem into three sub-problems where we optimize the three optimization variables separately. First, we assign the SF while respecting the LoRa specifications and assuring either fairness or unfairness between LoRa users. Second, we optimize the EH time using a one-dimensional exhaustive search. Then, we optimally optimize the transmit powers for the given SF and EH time. 
   
    \par}
\subsection{Unfair and Fair SF Allocation Algorithms}
    {\ First, we investigate the spreading factor allocation between the LoRa users. In order to assign the spreading factors between the users, we need to satisfy some conditions from LoRa specifications. The received signal at the gateway should exceed its sensitivity. The receiver sensitivity of the gateway depends on SF as shown in Table I in~\cite{canlorascale,SF_collision_distance}. In the literature, there are different ways to allocate the SFs between the users, as was mentioned in Section I.  
    In~\cite{canlorascale,SF_collision_distance}, the SFs are allocated according to the distance between the users and the gateway.  
    In~\cite{fair_alloc}, the fair collision probability $\prob(\text{SF}=f)= \frac{f/2^f}{ \sum_{i=7}^{12} i/2^i }$, for $f=7,\dots,12$, has been proposed to avoid the near-far problems.
    In~\cite{explora}, the algorithm EXPLoRa-SF has been proposed to equally divide the SFs between the users while respecting the received signal strength in (RSSI) values and relevant constraints. 
    \par}
    {\ Aligned with all these schemes, we propose to two different types of algorithms that assure either fairness or unfairness between the users. First, we assume that all users transmit with the maximum possible power $P_{n,max} = \min\left(P_t, \frac{1-d}{d} E_n\right)$, $\forall n$. 
    The unfair SF allocation equally divides the users to 6 groups. The fair SF allocation uses the fair collision probability in~\cite{fair_alloc}. More details are described in Algorithm \ref{algoSF}. Note that $sensi_{min}$  is the minimum required sensitivity at the gateway corresponding to $\text{SF}=12$ in Table 1 in~\cite{canlorascale}. 
        \begin{algorithm} 
       \SetAlgoLined
        \KwData{$P_t$, $d$, $E_n$, $U$, $g_n$, $sensi_{min}$}
         Initialize $U_a=0$\;
         \For{$n=1\rightarrow U$}{
         Compute $RSSI_n= P_{n,max} g_n$\;
         \lIf{$RSSI_n\geq sensi_{min}$}{ Increment $U_a$}}
         Order the $U_a$ users s.t. $RSSI$s are in a descending way\;
         Divide the $U_a$ ordered users into 6 groups of size $k_f$\;
          \For{$f=1\rightarrow 6$}{
          \leIf{Unfair}{ $k_f=\frac{U_a}{6}$}{$k_f= \frac{\frac{f+6}{2^{f+6}} }{\sum_{i=7}^12 \frac{i}{2^i}} U_a$}
          \lFor{$n=\sum_{j=1}^{f-1} k_{j}+1 \rightarrow \sum_{j=1}^{f} k_{j} $}{$SF_n=f+6$}}
        \Return ${\text{SF}}_n$\; 
        \caption{Unfair/Fair SF Allocation Algorithms} 
        \label{algoSF}
        \end{algorithm}
    \remark{If we do not respect the LoRa specifications, the optimal solution for the SF assignment between the users that maximizes the minimum rate involves an exhaustive search over a number of possibilities equal to $6^U$. The number of users $U$ is expected to be very large, therefore, the exhaustive search would be practically prohibitive. In order to reduce this complexity, we opt for the two low complex SF allocation algorithms explained in Algorithm \ref{algoSF} which both respect the LoRa specifications and align with what was done before in the literature.
    }
    \par}
\subsection{EH Time Allocation}
   {\ Since we are considering a "harvest-then-transmit" protocol, there is no data transmission during the moment that we harvest the energy. Also, the optimal EH time $\tau_{e,n}$ solution to \eqref{mainProb} should verify the two constraints C$_2$ and C$_3$. 
	Let us denote by $\delta_{n,max}^{(1)}=T_{off,n}$ and $\delta_{n,max}^{(2)}=\min\left(\frac{P_t}{E_n},\frac{1-d}{d}\right)T_{a,n}$.
	\par}
	\begin{theorem}
		{\ If the collision time in \eqref{cola} and \eqref{colb} is a non-monotonic function of the EH time, the optimal EH time $\tau_{e,n}$ is obtained by a one-dimensional search method in the interval $\left(0,\delta_{n,max}^{(1)}\right]$. 
		If the collision time is a monotonic nonincreasing function, the optimal EH time $\tau_{e,n}$ is equal to $\delta_{n,max}^{(1)}$.  
		If the collision time is a monotonic nondecreasing function or independent of the EH time, the optimal EH time is equal to $\delta_{n,max}^{(2)}$. 
		\par}
   \end{theorem}
   \begin{IEEEproof}
      {\ Since we have a maximum transmit power constraint in C$_1$, we represent two possible cases: either we are harvesting more than what we need (i.e. $P_t \leq \frac{\tau_{e,n}}{T_{a,n}} E_n$) or we are harvesting less than what we need (i.e. $P_t >\frac{\tau_{e,n}}{T_{a,n}} E_n$). 
		If $P_t \leq \frac{\tau_{e,n}}{  T_{a,n}}  E_n$, the constraint C$_2$ is trivially satisfied and the EH time should satisfy
		\begin{align}
			& \delta_{n,max}^{(2)} \leq \tau_{e,n} \leq \delta_{n,max}^{(1)} .\label{optEHsola}
		\end{align}
		At this point, one would say that it is wasteful harvesting more than what we need (i.e. $P_t$). However, spending more time on harvesting may reduce/increase the collision time between users and hence increase/decrease the throughput rate. Since our interest is minimizing the collision, optimizing the EH time depends on the monotonicity of the collision time in the EH time function. 
		If the collision time is a non-monotonic function with respect to the EH time, the optimal EH time is obtained by a one-dimensional search method in $\left[\delta_{n,max}^{(2)}, \delta_{n,max}^{(1)}\right]$, such as the bisection method \cite{cvx}. 
		
		If the collision time is a monotonic nonincreasing function with respect to the EH time, the optimal EH time is $\delta_{n,max}^{(1)}$. 
		If the collision time is a monotonic nondecreasing function with respect to the EH time (or independent of the EH time), the optimal EH time is $\delta_{n,max}^{(2)}$. 
		Otherwise, the optimal value is between $P_t \frac{T_{a,n} }{E_n}$ and $\frac{1-d}{d} T_{a,n}$. \\
		On the other hand, if $P_t >\frac{\tau_{e,n}}{  T_{a,n}}  E_n$, the EH time satisfies
		\begin{align}
			& 0 < \tau_{e,n} \leq \delta_{n,max}^{(2)}\leq \delta_{n,max}^{(1)}.\label{optEHsolb}
		\end{align}
		Similarly, the EH time can be obtained by a one-dimensional search method and the optimal value of the EH time depends on the collision time. 
		If the collision time is a nonincreasing function with respect to the EH time, the optimal EH time is $\delta_{n,max}^{(1)}$, which is the maximum value that we can consider. 
		If the collision time is a nondecreasing function with respect to the EH time (or independent of the EH time), the optimal EH time is $ \delta_{n,max}^{(2)}$.
        \par}
	\end{IEEEproof}
	\begin{corollary}
	    {\ If the collision time is given by the worst case scenario in \eqref{colWC}, the optimal solution of the EH time $\tau_{e,n}$ is exactly given by $\delta_{n,max}^{(2)}$, which is proportional to $T_{a,n}$. The available power after harvesting at each LoRa user is independent of $T_{a,n}$. Thus, the preference of SF has no effect on the energy harvesting constraint at each LoRa.
	    \par}
	\end{corollary}
	\begin{IEEEproof}
		{\ The proof is an immediate result of Theorem 1.
		\par}
	\end{IEEEproof}
\subsection{Optimal Power Allocation}
    {\ Given $\tau_{e,n}$ and $SF_n$, the optimal transmit powers $p_n$s at all LoRa users are solutions to
    	\begin{subequations}
        		\begin{align}
        		\underset{p_n}{\max} \hspace{3mm} \underset{n\in U_a}{ \min} \hspace{3mm} &  \log\left(1 + \gamma_n^{col} \right) ,\\
        		\text{s.t. }
        		& 0 \leq p_n \leq \overline{P}_n= \min\left(P_t,  \frac{\tau_{e,n}}{  T_{a,n}}  E_n\right).
        		\end{align}
        	\end{subequations}
    This problem is non convex therefore in order to solve this problem, we introduce a new optimization variable $t$ and we solve this equivalent convex problem:
        \begin{subequations}
        		\begin{align}
        		\underset{t, p_n, \forall n \in U_a}{\max } \hspace{3mm} & t,\\
        		\text{s.t. }
        		& 0 \leq p_n \leq \overline{P}_n,\\
		& t \leq  \log\left(1 + \gamma_n^{col} \right), \\
		& t\geq 0
        		\end{align}
        \end{subequations}	
	For a given $t_{low} \leq t \leq t_{up}$, the optimization problem is convex and the powers $p_n$'s can be optimally obtained for a given t. Then, the optimal $t$ can be selected using any one-dimensional search method such as the bisection method \cite{cvx}.
    \par}
\section{Numerical Results}
    {\ In this section, we present some selected simulation results to validate our proposed solution.
    \par}
    \begin{figure}[t]
         \begin{center}
               \includegraphics[scale=0.175]{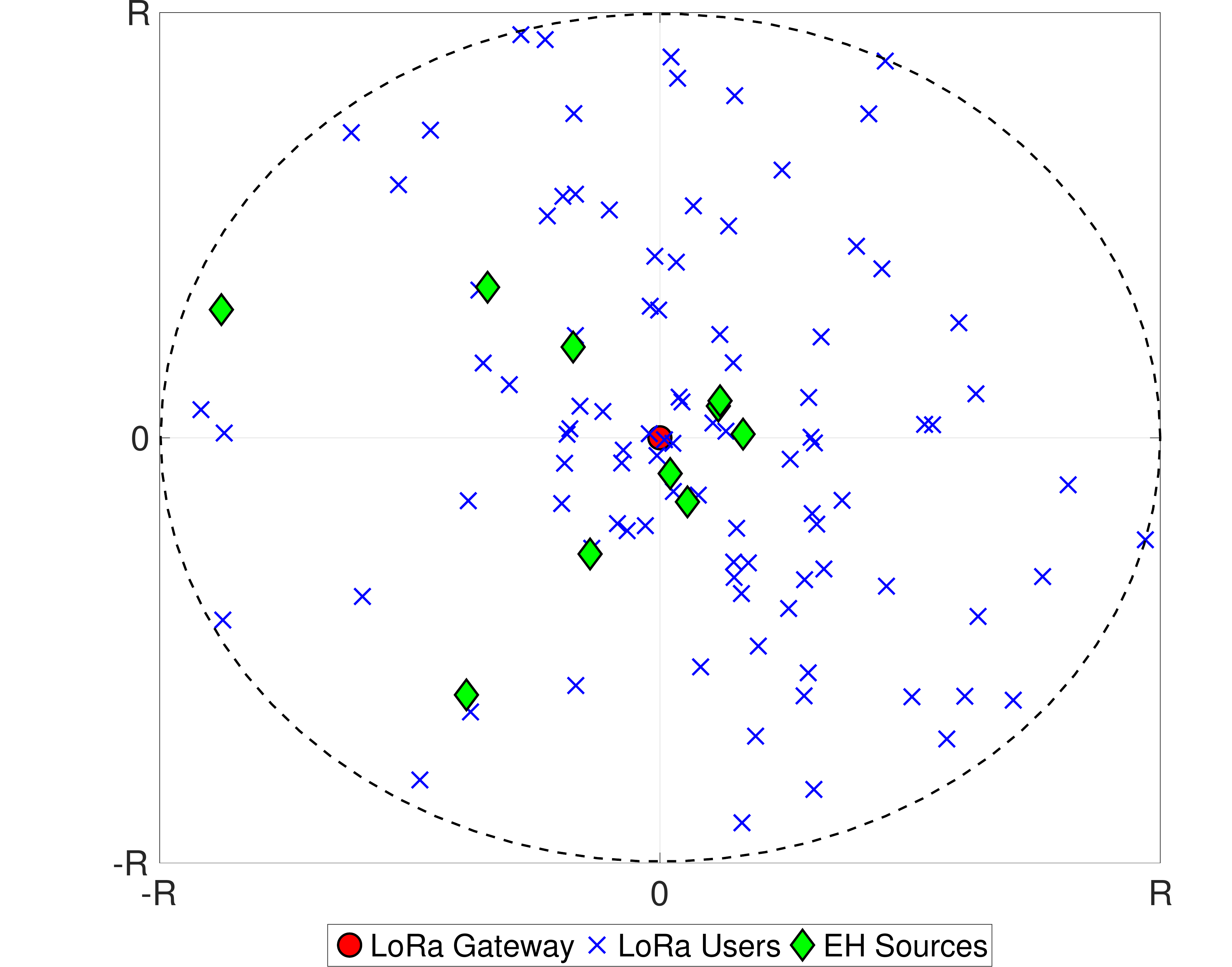}
               \caption{An example of LoRa network consisted of one gateway and a number of LoRa users with density $10^{4}$ Users/km$^2$ harvesting from power beacons with density $10^{3}$ PBs/km$^2$.}
    		   \label{network_LoRa_EH}
         \end{center}
    \end{figure} 
\subsection{LoRa Simulation Parameters}
    {\ For LoRa network, we assume that the users are uniformly distributed in a circle centred around the gateway. The simulation parameters are chosen following the LoRa specifications \cite{doppler}. 
    The noise variance is defined as $\sigma^2= -174+\text{NF} + 10\log_{10}(\text{BW})$ in dBm, where $\text{NF}$ is the noise figure equal to $6$ dB. The number of symbols $nb_n$ for each user is defined as $nb_n=n_{PR}+n_{PL,n}+4.25$, where $n_{PR}$ is the number of symbols in the preamble chosen equal to $12.25$, $n_{PL,n}$ is the number of symbols in the payload equal to $8+max\left(ceil\left(\frac{8\text{PL}-4\text{SF}_n+28+16}{ 4\left(\text{SF}_n-2\text{DE}\right)}\right)\left(\text{CR}+4\right),0\right)$, where the number of payload bytes is $\text{PL}=10$, the coding rate is $\text{CR}=1$, $\text{DE}=1$ for $\text{SF} \in\{11,12\}$ and $\text{DE}=0$ for $\text{SF}\in \{7,\dots,10\}$. 
    The path-loss exponent for both the power transfer and the information transfer links is $3.5$. The maximum transmit power for all LoRa users is $P_t=17$ dBm.
    \par}
\subsection{RF Energy Harvesting Model}
    {\ In our simulations, we are harvesting from the ambient RF signals transmitted by $N_b$ power beacons (PBs) randomly located in the cell with radius $R$. The transmit power of PBs is $P_b=1$ Watts. An example of the network is presented in Fig. \ref{network_LoRa_EH} where the power beacons density is $10^3$ PBs/km$^2$ and the LoRa users density is $10^4$ Users/km$^2$. Moreover, the harvested energy at each LoRa user depends on which EH model is considered either linear or nonlinear~\cite{fundamsurveryEH}. 
    In our analysis, we choose to follow the nonlinear EH model that has been proposed in~\cite{elena} where a sigmoidal model was shown to fit the experimental data. For the sigmoidal model, if $P_{rec,n}$ is the received power from the PBs at the user $n$, the harvested energy $E_n$ is expressed as
        \begin{align}
            E_n & = \Psi\left(P_{rec,n}\right), 
        \end{align}
    where $\Psi(\cdot)$ is the function defined as $\Psi(x)  = \frac{\beta(x) -M \Omega}{1-\Omega}$, $\beta(x)= \frac{M}{1+e^{-a(x-b)}}$, $\Omega=\frac{1}{1+e^{ab}}$, $M$ is the maximum harvested energy, $a$ and $b$ are experimental parameters which reflect the nonlinear charging rate with respect to the input power and the minimum required turn-on voltage for the start of current flow through the diode, respectively. In all plotted figures, the nonlinear model is considered, with the parameters $a=1500$, $b=0.0022$, and $M=24$ mW which were shown in \cite{elena} to fit the experimental data in \cite{datael}. 
    \par}
\subsection{SF Allocation Algorithms: Fairness or Unfairness?}
    \begin{figure}[t]
    	\begin{center}
            \includegraphics[scale=0.22]{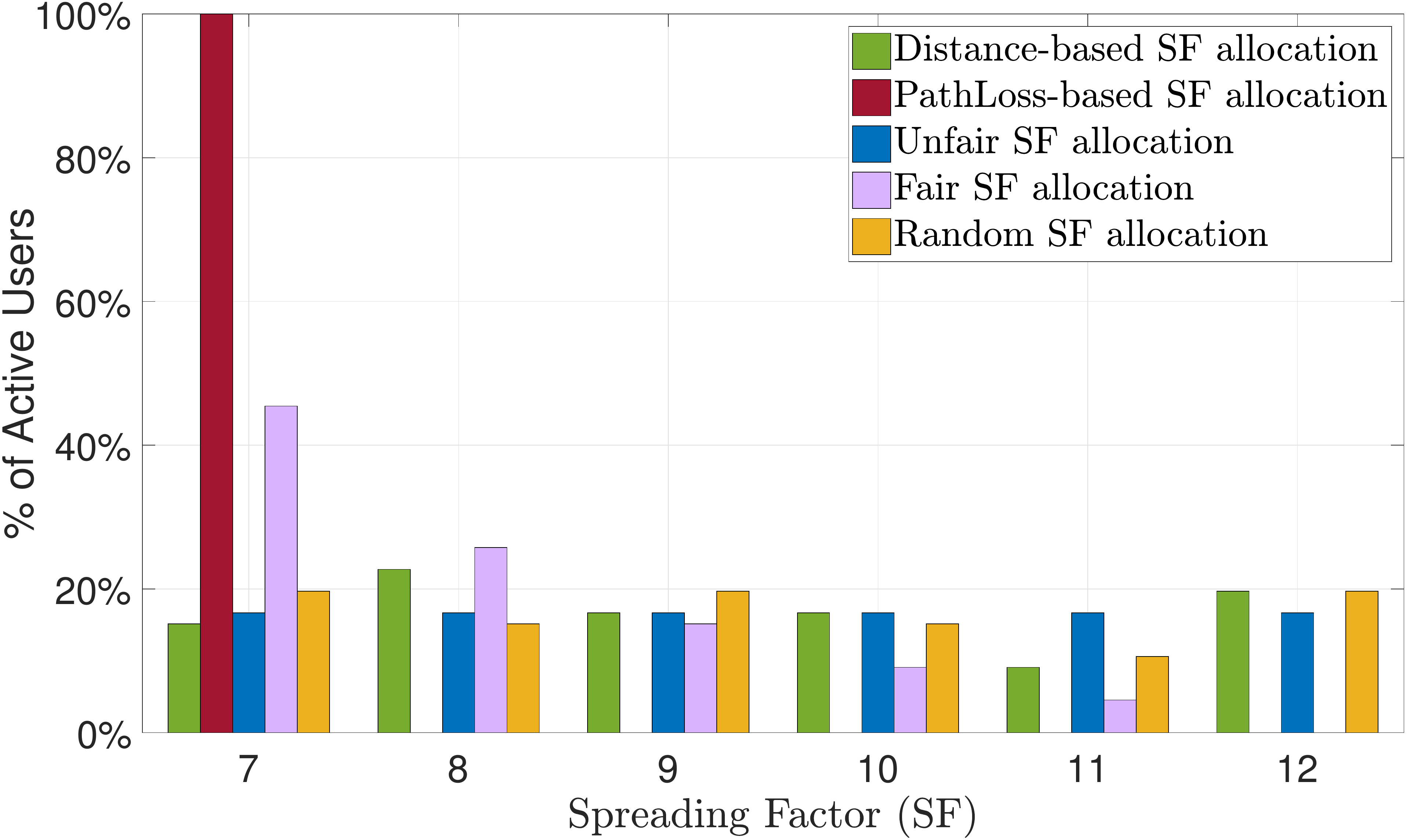}
          	\caption{Allocation of active users to different SFs, $SF\in\{7,\dots,12\}$, depending on the different SF algorithms, for a total of active users with density $10^{4}$ Users/km$^2$.}
         	\label{users_vs_SFalloc}
       \end{center}
    \end{figure} 
    {\ In Fig. \ref{users_vs_SFalloc}, we show the assignment of active users to different SFs, $\text{SF}\in\{7,\dots,12\}$, depending on the different SF allocation algorithms. The fair and unfair SF allocation algorithms are the ones described in Algorithm \ref{algoSF}. 
    The distance-based SF allocation assigns the SFs depending on the distance between the users and the gateway~\cite{canlorascale}. Let $d_i=\frac{iR}{6}$, for $i=0,\dots,6$. The users between the distances $d_i$ and $d_{i+1}$ will be assigned $\text{SF}=7+i$, for $i=0,\dots,5$. The pathloss-based SF allocation assigns the users while respecting Table I in~\cite{canlorascale,SF_collision_distance}. 
    We can observe that the pathloss-based SF allocation, and the fair SF allocation tend to assign more users to lower SFs, while the other algorithms assign the users to all SFs. We will see in the following figures either the fairness or unfairness is better for the minimum rate performance of the LoRa users.
    \par}  
\subsection{Optimal vs Maximum EH Time}
    \begin{figure}[t]
        \begin{center}
          	\includegraphics[scale=0.23]{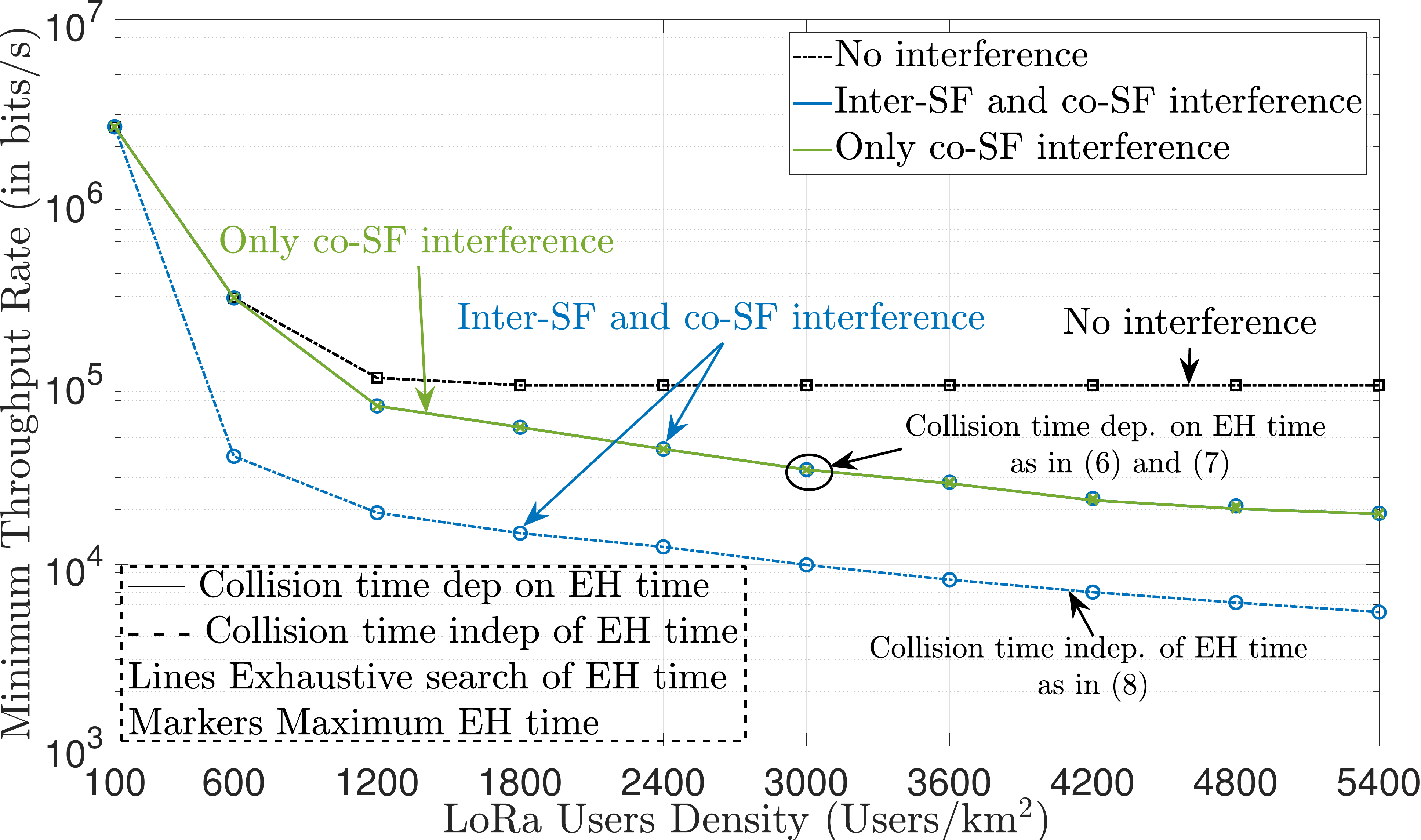}
          	\caption{The minimum throughput rate of LoRa users versus the density of active LoRa users per km$^2$ with unfair SF allocation algorithm for different multiuser interference scenarios, and the packet collision time either depending or not on EH time, with $P_t= 17$ dBm.}
         	\label{R_U_exhaus}
       \end{center}
    \end{figure} 
    {\ Before comparing the SF allocation algorithms and the multiuser interference scenarios, we study first the choice of the EH time. 
    In Fig. \ref{R_U_exhaus}, we have plotted the minimum throughput rate of LoRa users versus the density of active LoRa users per km$^2$ with unfair SF allocation algorithm. We consider different multiuser interference scenarios. 
    The no interference case refers to the case where the correlation factor $\rho_{n,m}=0$, $\forall n, m$. The only co-SF interference case refers to the case where the correlation factor $\rho_{n,m}=1$ if the users $n$ and $m$ have the same spreading factor and where the correlation factor $\rho_{n,m}=0$ if the users $n$ and $m$ have different spreading factors. 
	The inter-SF and co-SF interference case refers to the worst scenario case where $\rho_{n,m}=1$, $\forall n, m$. 
    The packet collision time is either depending on EH time as in \eqref{cola} and \eqref{colb} or not depending on EH time as in \eqref{colWC}. 
    We compare the results using the optimal EH time obtained using the exhaustive search versus the maximum EH time. The maximum value of EH time is equal to $\delta_{n,max}^{(1)}$ for the collision time depending on EH time; and equal to $\delta_{n,max}^{(2)}$ for the collision time not depending on EH time in \eqref{colWC}. 
    The lined curves refer to the exhaustive search solution of EH time and the marked curves refer to the maximum value of EH time. In all plotted scenarios, we can see the agreement between the exhaustive search and the maximum value of EH time. Thus, in the following figures, we will consider the special case where the EH time is equal to its maximum value. 
    \par}
\subsection{Collision Time and Multi-User Interference}
    \begin{figure}[t]
       \begin{center}
          	\includegraphics[scale=0.23]{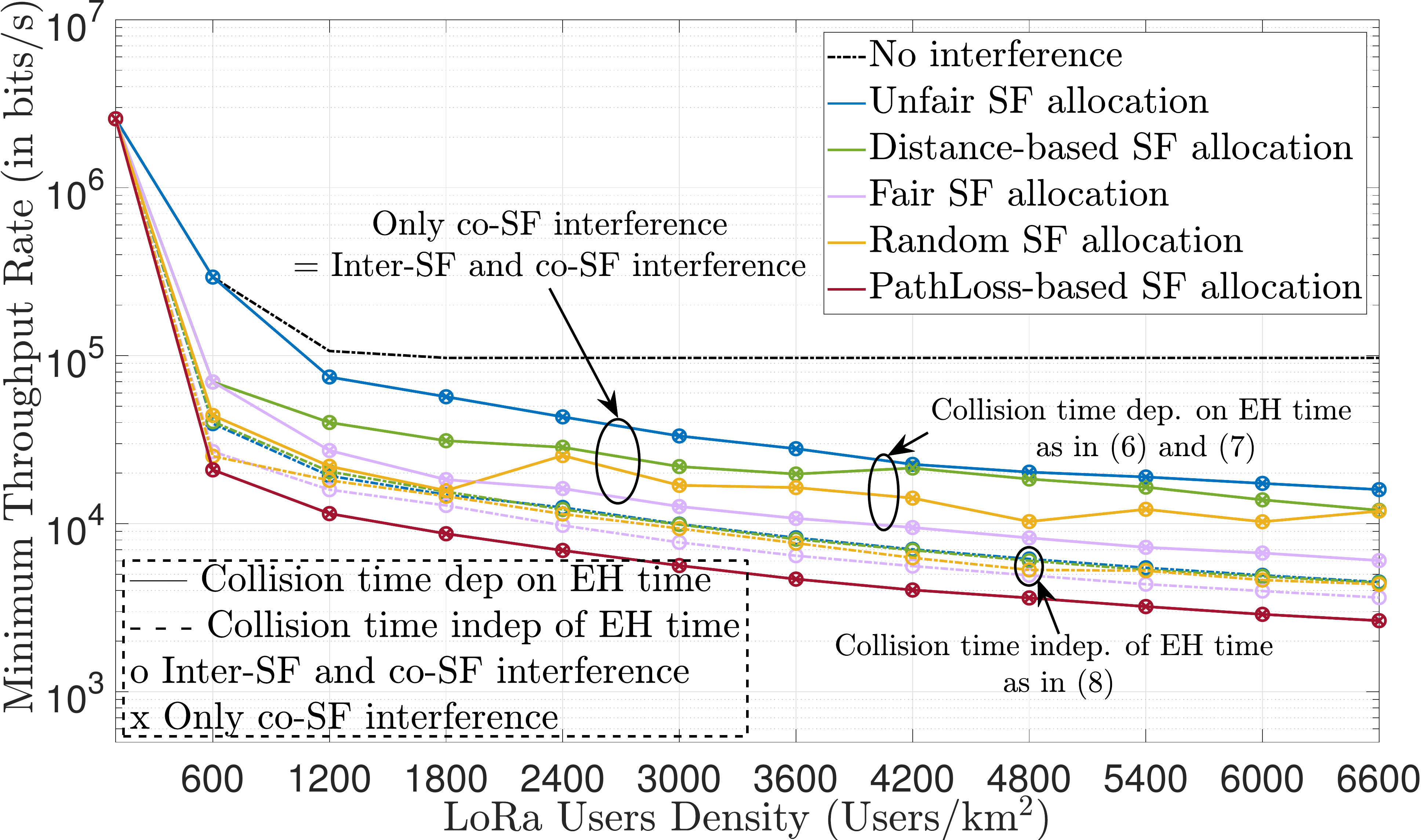}
          	\caption{The minimum throughput rate of LoRa users versus the density of active LoRa users per km$^2$, for different SF allocation algorithms, different multiuser interference scenarios, and the packet collision time either depending or not on EH time, with $P_t= 17$ dBm.}
         	\label{R_nbr}
      \end{center}
    \end{figure} 
	{\ In Fig. \ref{R_nbr}, we have plotted the minimum throughput rate of LoRa users versus the density of active LoRa users per km$^2$ for different SF allocation algorithms and different multiuser interference scenarios. The packet collision time is either depending or not on EH time. 
	The lined curves are obtained with a packet collision time depending on the EH time as in \eqref{cola} and \eqref{colb}. 
	The dashed curves refer to the packet collision time independent of the EH time as in \eqref{colWC}. 
	First of all, we can see that when the collision time is defined as in \eqref{cola} and \eqref{colb}, the minimum rate outperforms the scenario where the collision time is defined as in \eqref{colWC}. This observation is expected since fewer users are transmitting simultaneously when the collision time depends on the EH time. Whereas, all users are transmitting simultaneously in \eqref{colWC}. 
	In addition, note that when $\tau_{e,n}=\delta_{n,max}^{(1)}$, the collision time simplifies to either zero if the users have different spreading factors, and to the time on air if the users have the same spreading factor, as shown in Remark \ref{rk6}. This means that, for this case, we don't have an inter-SF interference. That is why we have an agreement between the minimum rate with the inter-SF and co-SF interference and the one with only co-SF interference, as shown in Fig. \ref{R_nbr}. 
	Hence, we can conclude that the expression of the collision time depending on the EH time protects the system from the inter-SF interference even though we consider imperfect orthogonality between the SFs. The remaining limitation of performance here is only due to the co-SF interference. 
	\par}
	{\ Furthermore, the unfair SF allocation with the collision time expressed in \eqref{cola} and $\tau_{e,n}=\delta_{n,max}^{(1)}$ reaches the nearest performance to the no interference case and outperforms the other SF allocation algorithms. 
	The pathloss-based SF allocation has the poorest performance in terms of minimum rate. As shown in Fig. \ref{users_vs_SFalloc}, this algorithm tends to assign more users to the lower SFs. Thus, we can conclude that the best performance is achieved while equally allocating the SFs between users and not while assigning the nearest users to the lowest SFs, in contrast with what was shown previously in \cite{fair_alloc}. 
\par}

\section{Conclusion}
    {\ In this paper, we have studied the uplink resource allocation in LoRa networks powered by ambient energy harvesting. First, the packet collision time between users using the same or different SFs was expressed in function of the EH time duration. We proposed two types of SF allocation algorithms which rely on either the fairness or unfairness between the users. Through the simulation results, we have seen that using the unfair SF allocation and the collision time depending on the EH time outperforms all the other scenarios. We have seen also that the expression of the collision time cancelled the effect of inter-SF interference when the EH time is equal to the off time. 
    Finally, we have seen that the co-SF interference is the main limitation of the throughput performance, not really the energy scarcity. 
    We can conclude on the importance of the transmission scheduling on the number of colliding users and hence the performance of the minimum throughput rate. Also, it seems appealing to consider the multiuser interference management techniques to cancel the co-SF interference. 
   
    \par}
\bibliographystyle{IEEEtran}
\bibliography{references}
\end{document}